\documentclass{article}                                                                   

\usepackage{amsfonts, graphicx}
\usepackage{latexsym}                                                                      
\usepackage{amsmath,amsthm,amsfonts}                                              
\usepackage{amssymb}                                                                        
\usepackage{palatino}                                                                       
\usepackage[mathscr]{eucal}

\newcommand{\R}{\Bbb R}

\def\qed{\hbox to 0pt{}\hfill$\rlap{$\sqcap$}\sqcup$}
\newtheorem{theorem}{Theorem}[section]

\theoremstyle{definition}

\theoremstyle{remark}

\newcommand{\inserttitle}[1]{\begin{center}{\Large\bf #1}\end{center}} 
\newcommand{\insertauthor}[1]{\begin{center}{\sc #1}\end{center}}
\numberwithin{equation}{section}
\begin{document}

\inserttitle{On the minimal speed of traveling waves for a non-local
delayed reaction-diffusion equation}

\insertauthor{Maitere Aguerrea and Gabriel Valenzuela}

{\footnotesize \centerline{Instituto de Matem\'atica y F\'{\i}sica}
\centerline{ Universidad de Talca, Casilla 747, Talca, Chile } }

\vspace{0.5 cm}



\bigskip
\begin{quote}{\normalfont\fontsize{8}{10}\selectfont
{\bfseries Abstract.} In this note, we give constructive upper and
lower bounds for the minimal speed of propagation of traveling waves
for non-local delayed reaction-diffusion equation.
\par}
\end{quote}

\section{Introduction and the main results}
In this note, we estimate
the minimal speed of propagation of positive traveling wave
solutions for non-local delayed reaction-diffusion equation
\begin{equation}\label{17}
u_t(t,x) = u_{xx}(t,x)  - u(t,x) + \int_{\R}K(x-s)g(u(t-h,s))ds, \ u
\geq 0,\ x \in \R,
\end{equation}
which is widely used in applications, e.g. see \cite{gouss,LRW,SZ,
STPA, WLR} and references wherein. It is assumed that the birth
function $g$ is of the monostable type, $p: =g'(0) >1$ and $h\geq0$.
The non-negative kernel  $K$ is such that $K(s)=K(-s)$ for $s\in \R$,
$\int_{\R}K(s)ds=1$ and $\int_{\R}K(s)\exp(\lambda s)ds$ is finite
for all $\lambda \in \R$. Consider
\begin{equation}\label{ce}
\psi(z,\varepsilon)=\varepsilon
z^2-z-1+p\exp(-zh)\int_{\R}K(s)\exp(-\sqrt{\varepsilon}zs)ds,
\end{equation}
which determines the eigenvalues of  Eq. (\ref{17}) at the trivial
steady state. From \cite{M,STPA}, we know that there is
$\varepsilon_{0}= \varepsilon_{0}(h)
> 0$ such that $\psi(z,\varepsilon_{0})=0$ has a unique multiple
positive root $z_{0}= z_{0}(h)$. Furthermore, if $g(s)\leq g'(0)s$ for $s\geq0$, then the minimal speed $c_{*}$ is equal to $c_{*}=1/\sqrt{\varepsilon_{0}}$ . Note that $z_{0}$ and
$\varepsilon_{0}$ are the unique solutions of the system
\begin{equation}\label{se}
\psi(z,\varepsilon)=0, \quad  \psi_{z}(z,\varepsilon)=0.
\end{equation} It is known \cite{WLR} that for
various systems modeled by equation (\ref{17}), the minimal wave
speed $c_*$ coincides with the spreading speed. Therefore, it is important to study the effects caused by the delay and other
parameters (depending on specific models) on $c_*$,  cf.
\cite{LRW,S,SZ,ETST,WLR}. Another aspect of the problem
concerns easily calculable upper and lower bounds for $c_*$.  In
particular, in the recent work \cite{WM}, Wu \textit{et al.} give several nice estimations for $c_*$ when
$K_{\alpha}(s)=\frac{1}{\sqrt{4\pi \alpha}}e^{-s^2/4\alpha}$ and
$\alpha \leq h$. However, the
approach of \cite{WM} depends heavily on the condition $\alpha\leq h$ and on the special form of $K$ which is the
fundamental solution of the heat equation. In the present work, we
use a completely different idea to estimate the minimal speed for
general kernels and without imposing any restriction on $h$.

Let us state our main result. Set

\begin{eqnarray*}
k_{1}&=&2\sqrt{\frac{p-1}{1+\frac{p}{2}\int_{\R}\ s^2K(s)ds}}-p\int_{\R}sK(s) \exp\left(-s\sqrt{\frac{p-1}{1+\frac{p}{2}\int_{\R}\ s^2K(s)ds}}\right)ds,\\
k_{2}&=&\frac{1}{\sqrt{\ln p}}\ln{\left( p\int_{\R}K(s)\exp(-\sqrt{\ln p}s)ds \right) }.
\end{eqnarray*}
It is clear that $k_{2}>0$ and bellow we will show that $k_{1}$ is positive.
\begin{theorem}\label{th} Assume that $K(s)\geq0$ is such that $K(s)=K(-s)$ for $s \in \R$,
$\int_{\R}K(s)ds=1$ and $\int_{\R}K(s)\exp(\lambda s)ds$ is finite
for all $\lambda \in \R$. Then $c_{*}=c_{*}(h)=1/\sqrt{\varepsilon_{0}(h)}$ is a $C^{\infty}$-smooth decreasing function of variable $h \in \R_{+}$. Moreover, 
\begin{enumerate}
\item $ \displaystyle \max \left\{ 2\sqrt{ \frac{{p-1}}{p(2h+h^2)+1}}, \frac{2\sqrt{\ln{p}}}{1+h} \right\}< c_{*}<\min \left\{ \frac{k_{1}}{1+h}, \frac{k_{2}}{h} \right\}$, $h \in [0,1]$,\\

\item $ \displaystyle \max \left\{ 2\sqrt{ \frac{{p-1}}{p(2h+h^2)+1}}, \frac{\sqrt{\ln{p}}}{h} \right\}< c_{*}<\min \left\{ \frac{k_{1}}{2}, \frac{k_{2}}{\sqrt{h}} \right\}$, $h \in [1,+\infty)$.
\end{enumerate}

Furthermore, $\displaystyle  \frac{C_{1}}{h}\leq c_{*}(h) \leq \frac{C_{2}}{h}$, \ $h\geq1$, for some positive $C_{1}<C_{2}$.
\end{theorem}

Observe that Theorem \ref{th} implies that $c_{*}(h)=O(h^{-1})$, $h \to +\infty $, in this way we improve the estimation $c_{*}(h)=O(h^{-1/2})$, $h \to +\infty $, proved in \cite{ETST,WM}.

\begin{proof} It follows from \cite{M, STPA} that the functions $z_{0}= z_{0}(h)$ and $\varepsilon_{0}= \varepsilon_{0}(h)$ are well defined for all $h\geq0$. Set $F(h,z,\varepsilon)=(\psi(z,\varepsilon),
\psi_{z}(z,\varepsilon))$. It is easy to see  $F\in C^{\infty}(\R_{+}\times \R
\times (0, \infty), \R^2), \ F(h,z_0,\varepsilon_0)=0$, and
\begin{eqnarray*}
\left|\frac{\partial F(h, z_0,\varepsilon_0)}{\partial(
z_0,\varepsilon_0)}\right|&=&
\psi_{zz}(z_{0},\varepsilon_{0})\psi_{\varepsilon}(z_{0},\varepsilon_{0}) \\
&=&(2\varepsilon_{0}+p\int_{\R}K(s)\exp({-z_{0}(h+\sqrt{\varepsilon_{0}}s))(h+\sqrt{\varepsilon_{0}}s)^2}ds) \\
&\times& \frac{z_{0}}{2\varepsilon_{0}}(1+hp\int_{\R}K(s)\exp(-z_{0}(h+\sqrt{\varepsilon_{0}}s))ds) >0.
\end{eqnarray*}
Applying the Implicit Function Theorem we find that $z_{0}, \varepsilon_{0} \in C^{\infty}(0,+\infty)$.

On the other hand, after introducing a new variable $w=\sqrt{\varepsilon}z$ we find that system (\ref{se}) takes the following form:

\begin{equation}\label{ew}
\left(1+\frac{w}{\sqrt{\varepsilon}}-w^2\right)\exp\left(\frac{wh}{\sqrt{\varepsilon}}\right)=p\int_{\R}K(s)\exp(-ws)ds.
\end{equation}
\begin{equation}\label{eww}
\left(\frac{h}{\sqrt{\varepsilon}}w^2+\left(2-\frac{h}{\varepsilon}\right)w-\frac{1+h}{\sqrt{\varepsilon}}\right)\exp\left(\frac{wh}{\sqrt{\varepsilon}}\right)=p\int_{\R}sK(s)\exp(-ws)ds.
\end{equation}
Let $G(w)=\left(1+\frac{w}{\sqrt{\varepsilon_{0}}}-w^2\right)$, $H(w)=\left(1+\frac{w}{\sqrt{\varepsilon_{0}}}-w^2\right)\exp\left(\frac{wh}{\sqrt{\varepsilon_{0}}}\right)$ and $R(w)=p\int_{\R}K(s)\exp(-ws)ds$. Set also $w_{0}=w_{0}(h)=\sqrt{\varepsilon_{0}(h)}z_{0}(h)$. First, note that $G(w_{0})=\exp\left(\frac{-wh}{\sqrt{\varepsilon_{0}}}\right)R(w_{0})>0$ and $G(w)\geq 1$ when $0\leq w\leq 1/\sqrt{\varepsilon_{0}}$. As can be checked directly, $H$ has a unique positive local extremum (maximum) at some $\bar{w}$. Since $K(s)=K(-s), \ s \in \R$, it is easy to see that $R$ increases on $\R_{+}$.

Differentiating equation (\ref{ew}) with respect to $h$ and using (\ref{eww}) we get the following differential equation

\begin{equation}\label{eh}
\varepsilon'_{0}(h)=\frac{2 \varepsilon_{0}(h) G(w_{0}(h))}{1+hG(w_{0}(h))}>0.
\end{equation}

The remainder of the proof will be divided in
several steps.

\underline{Step I.} If $h \in [0,1]$, then $H'(1/\sqrt{\varepsilon_{0}})=\left(\frac{h-1}{\sqrt{\varepsilon_{0}}}\right)e^{h/\varepsilon_{0}}\leq0$. Hence, $\bar{w}\leq
1/\sqrt{\varepsilon_{0}}$. In addition,  if $w \in (0,\bar{w})$ then $H'(w)>0$. As $R'(w)>0$ for $w>0$, we have $w_{0}<\bar{w}\leq
1/\sqrt{\varepsilon_{0}}$. Thus, we get $G(w_{0})\geq 1$. In this way, $\varepsilon'_{0}(h)\geq
2\varepsilon_{0}(h)/(1+h)$ for $h \in [0,1]$ that yields
$(1+h)^2\varepsilon_{0}(0) \leq \varepsilon_{0}(h)\leq (1+h)^2\varepsilon_{0}(1)/4$ (equivalently,
$2c_*(1)/(1+h)\leq c_*(h)\leq c_*(0)/(1+h)$, for $ h \in [0,1]$).  Next, taking $h=0$
in equations (\ref{ew}) and (\ref{eww}) we obtain that
\begin{equation}\label{ewww}
\frac{1}{\sqrt{\varepsilon_{0}(0)}} = 2w_0(0)-
p\int_{\R}sK(s)\exp{(-w_0(0)s)}ds,
\end{equation}
\begin{eqnarray*}
1+w_{0}^2(0)&=&p\int_{\R}K(s)(1+w_{0}(0)s)\exp(-w_{0}(0)s)ds \\
&=& p\left(1-\frac{\int_{\R}s^2K(s)ds}{2} w_{0}^2(0)- \frac{\int_{\R}s^4K(s)ds}{8} w_{0}^4(0)- \dots\right).
\end{eqnarray*}
As a consequence of the latter formula, we get
$$
w_{0}(0)<\sqrt{\frac{p-1}{1+\frac{p}{2}\int_{\R}\ s^2K(s)ds}} \ .
$$
Then  (\ref{ewww}) implies that
$c_{*}(0)<k_1$ so that $c_{*}(h)<
k_{1}/(1+h)$ for $h\leq 1$. Note that $k_{1}>0$ since $R$ is increasing for $w>0$. Finally, since $c_{*}(h)$ is
decreasing, we have that $c_{*}(h)< k_{1}/2$ for $h\geq1$.

\underline{Step II.} 
If $h\geq1$, then $\bar{w}\geq 1/\sqrt{\varepsilon_{0}}$. As consequence, $G(\bar{w})\leq1=G(1/\sqrt{\varepsilon_{0}})$ so that $G(w)\geq G(\bar{w})$ for all $w \in [0,\bar{w}]$ (see Figure \ref{f1}). Additionally, $G(\bar{w})=(2\bar{w}\sqrt{\varepsilon_{0}}-1)\frac{1}{h}\geq \frac{1}{h}$, therefore we conclude that $G(w_{0})\geq 1/h$. Hence, we have $\varepsilon_{0}'(h)\geq \varepsilon_{0}(h)/h$, so that $\varepsilon(h)\geq\varepsilon(1)h$ (equivalently,
$c_*(h)\leq c_*(1)/\sqrt{h}$) for $h\geq1$. Now, if $h=1$ we have $\bar{w}=1/\sqrt{\varepsilon_{0}(1)}$. Thus, taking $h=1$ and $w=\bar{w}$ in (\ref{ew}) we get $\exp(1/\varepsilon_{0}(1))= R(\bar{w})>R(0)=p$ that yields $\sqrt{\ln p}<1/\sqrt{\varepsilon_{0}(1)}=c_{*}(1)$.
On the other hand, for all $0\leq w<1/\sqrt{\varepsilon_{0}}$, we have
\begin{equation}\label{ie}
\exp\left(\frac{wh}{\sqrt{\varepsilon_{0}}}\right) < \left(1+\frac{w}{\sqrt{\varepsilon_{0}}}-w^2\right)\exp\left(\frac{wh}{\sqrt{\varepsilon_{0}}}\right) \leq p\int_{\R}K(s)\exp(-ws)ds.
\end{equation}
In particular, taking $h=1$ and $w=\sqrt{\ln p}$ in (\ref{ie}) we conclude that $c_{*}(1) < k_{2}$ so that $c_{*}(h)<k_{2}/\sqrt{h}$ for $h\geq 1$. Additionally, using $c_*(h)\geq 2c_*(1)/(1+h)$ obtained in step I, we also concluded that $c_*(h)>2\sqrt{\ln p}/(1+h)$, for $h\in [0,1]$.

\begin{figure}[h]
\begin{center}
\includegraphics[width=8cm]{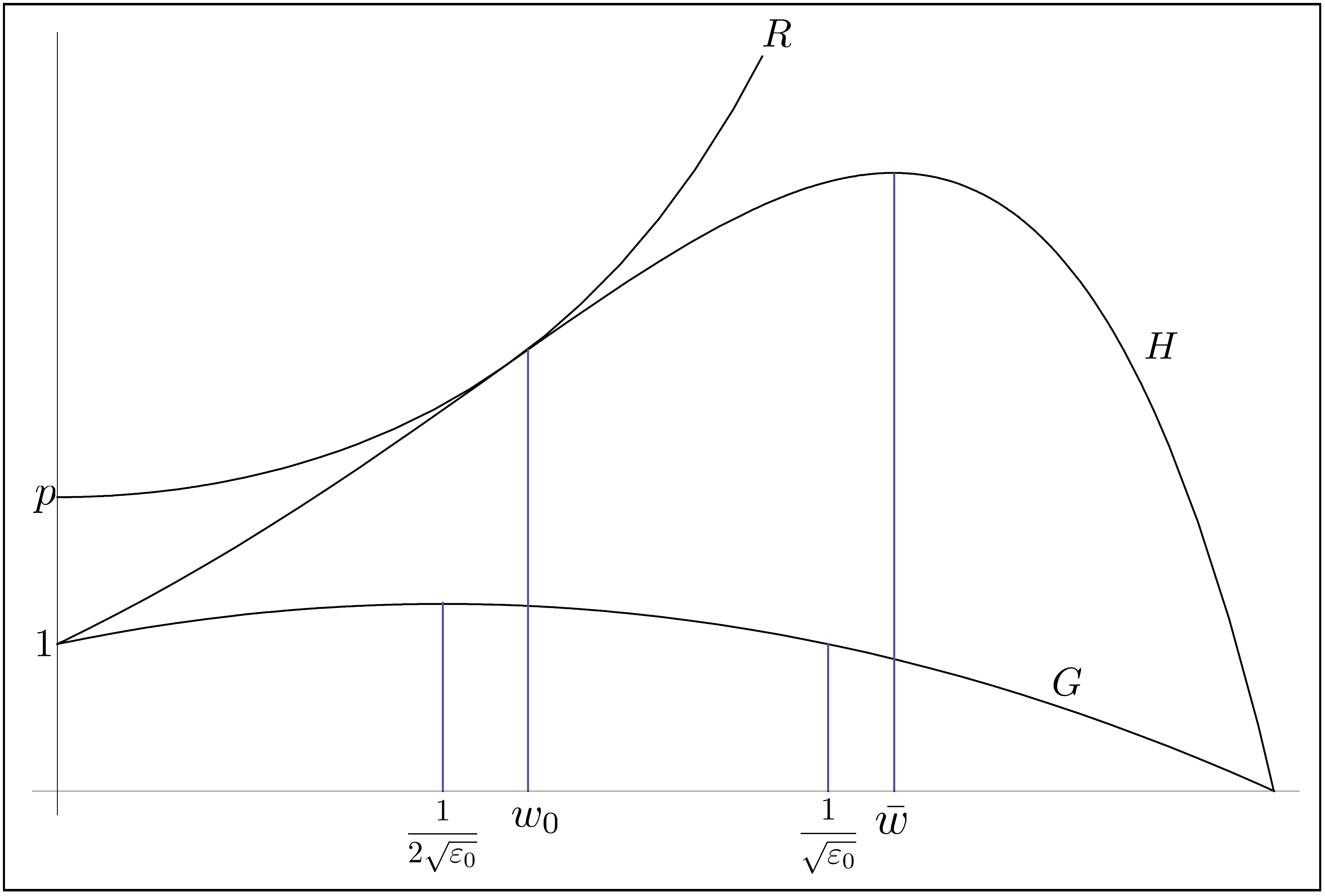}
\end{center}
\caption{$G$, $H$ and $R$ for $h>1$. }\label{f1}
\end{figure}

\underline{Step III.} For $h>0$, it is evident that
$\varepsilon_{0}'(h)\leq 2\varepsilon_{0}(h)/h$. Integrating the latter inequality on $[h,1]$ we obtain $\varepsilon(h)\geq\varepsilon(1)h^2$ (equivalently,
$c_*(h)\leq c_*(1)/h$), for $0<h\leq1$ so that $c_*(h)< k_{2}/h$, for $h \in (0,1]$.
Analogous, by integrating $\varepsilon_{0}'(h)\leq 2\varepsilon_{0}(h)/h$ on $[1,h]$ we have
$\varepsilon_{0}(h)\leq{\varepsilon_{0}(1)}h^2$ (equivalently,
$c_*(h)\geq c_*(1)/h$), for $h\geq1$. Thus, we obtain $c_{*}(h)>\sqrt{\ln{p}}/{h}, \ h\geq1$. 

On the other hand,  for all $h\geq0$, we have $G(w_{0})\leq 1+1/(4\varepsilon_{0})$. As consequence, $\varepsilon_{0}'(h)\leq(4 \varepsilon_{0}(h)+1)/(2(1+h))$  for all $h\geq0$ so that $\varepsilon_{0}(h)\leq((4\varepsilon_{0}(0)+1)(1+h)^2-1)/4$. Taking $h=0$ in (\ref{ew}), we get $1+1/(4\varepsilon_{0}(0))>G(w_{0}(0))=R(w_{0}(0))>p$ so that $c_{*}(0)>2\sqrt{p-1}$. In consequence,
\begin{equation}\label{add}
c_{*}(h)>2\sqrt{ \frac{{p-1}}{p(2h+h^2)+1}}, \quad h\geq0.
\end{equation}

\underline{Step IV.} Setting $w=r$, $r \in (0,1)$, in the second inequality of (\ref{ie}) we obtain
$$\left(1-r^2\right)\exp\left(\frac{rh}{\sqrt{\varepsilon_{0}(h)}}\right)< p\int_{\R}K(s)\exp(-rs)ds,$$
from which we get that
\begin{equation}\label{ad}
\frac{1}{\sqrt{\varepsilon_{0}(h)}}<\frac{1}{hr}\ln{\left( \frac{ p}{1-r^2}\int_{\R}K(s)\exp(-rs)ds\right)}, \quad h>0.
\end{equation}
Considering (\ref{add}) and (\ref{ad}) we get $\displaystyle \frac{C_{1}}{h}\leq c_{*}(h) \leq \frac{C_{2}}{h}$ for $h\geq1$. This completes the proof.
\end{proof}

\section{An example} 
Consider the heat kernel $K_{\alpha}(s)=(4\pi \alpha)^{-1/2}\exp{(-s^2/(4\alpha))}$. Then Theorem \ref{th} applies with 
\begin{eqnarray*}
k_{1}=2\sqrt{p-1} \left( \frac{1+\alpha p\exp{ \left( \frac{\alpha(p-1)}{1+\alpha p} \right) }}{\sqrt{1+\alpha p}} \right), \quad 
k_{2}= (1+\alpha)\sqrt{\ln{p}}.
\end{eqnarray*}
\begin{figure}[h]
\begin{center}
\includegraphics[width=9.6cm]{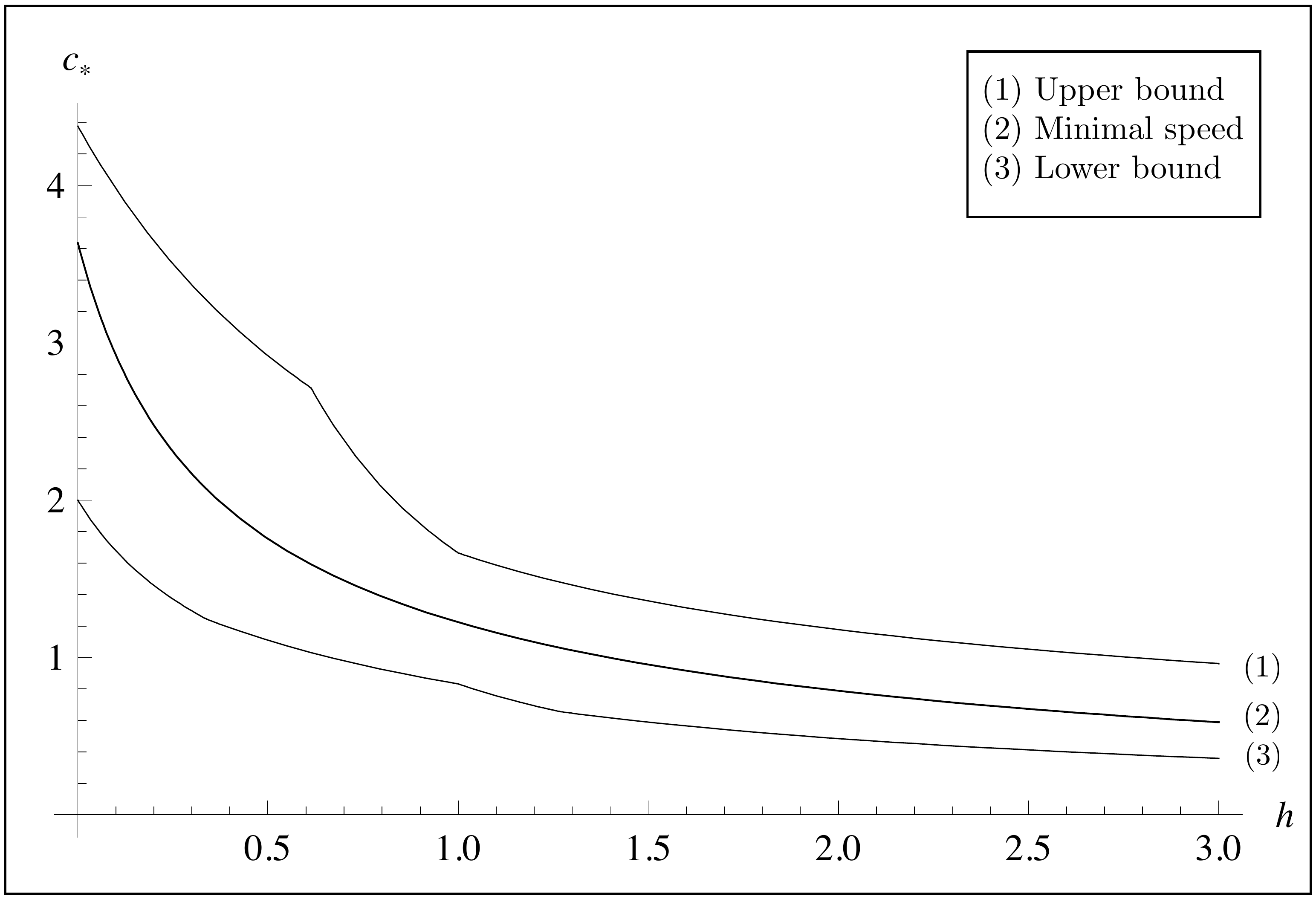}
\end{center}
\caption{The minimal speed and its bounds ($p=2$ and $\alpha=1$). }\label{f2}
\end{figure}
In fact, in this case we can plot graphs of $c_{*}$ against $h$ using standard numerical methods to solve some appropriately chosen initial value problem $\varepsilon_{0}(h_{0})=\rho_{0}$ for differential equation (\ref{eh}). For example, if we take $h_{0}=\alpha$ then $\rho_{0}$ coincides with positive solution of the equation
$1+ \frac{1}{4\rho}=p\exp{( - \frac{\alpha}{4\rho})}$. Next, we can explicitly find $G(w_{0})$ in (\ref{eh}) by using Cardano's formulas to solve the cubic equation
$(w_{0}^2-w_{0}/\sqrt{\varepsilon_{0}}-1)(2\sqrt{\varepsilon_{0}}\alpha w_{0}-h)+1-2\sqrt{\varepsilon_{0}} w_{0}=0.$
It is easy to see that this equation has three real roots for all $h\geq 0$ and $\alpha>0$, and that $w_{0}$ is the leftmost positive root.
 
Figure \ref{f2} shows the minimal speed $c_{*}$ and its estimations when $p=2$ and $\alpha=1$. Remark that we do not need the restriction  $\alpha\leq h$ required in \cite{WM}.

Finally, note that letting $\alpha \to 0^+$ in (\ref{17}) and (\ref{ce}) we recover the characteristic equation for the delayed reaction-diffusion equation
$$
u_{t}(t,x)=u_{xx}(t,x)-u(t,x)+g(u(t-h,x)),
$$
which was studied by various authors (e.g. see \cite{ATV, ETST} and references therein). In this case, our results complete and partially improve the estimations of \cite{ETST}.

\section*{Acknowledgments}
The authors thank Professor Sergei Trofimchuk for valuable discutions and helpful comments. This work was partially supported by  FONDECYT (Chile), project 1071053.


\begin{thebibliography}{99}

\bibitem{ATV} {M. Aguerrea, S. Trofimchuk and G. Valenzuela}, Uniqueness of fast travelling 
fronts in a single species reaction-diffusion equation with delay, Proc. R. Soc. 
A 464, published online: 20 May 2008.

\bibitem{gouss} {S. A. Gourley, J. So, J. Wu,} Non-locality of reaction-diffusion
equations induced by delay: biological modeling and nonlinear
dynamics, {J. Math. Sciences} 124 (2004) 5119-5153.

\bibitem{LRW} W.T. Li, S. Ruan, Z.C. Wang, On the diffusive Nicolson's Blowflies equation with nonlocal delays, {
J. Nonlinear Sci.,} 17 (2007) 505-525.

\bibitem{M} S. Ma, Traveling waves for non-local delayed diffusion equation via auxiliary equation, {
J. Differential Equations,}  237 (2007) 259-277.

\bibitem{S} K. W. Schaaf, Asymptotic behavior and traveling wave solutions for parabolic functional differential equations, {
Trans. Am. Math. Soc.,} { 302} (1987) 587-615.

\bibitem{SZ}  J. So, J.Wu and X. Zou, A reaction-diffusion model for a single species with age structure. I. Travelling wavefronts on unbounded domains, Proc. R. Soc. A 
457 (2001) 1841-1853.

\bibitem{STPA} {  E. Trofimchuk, P. Alvarado, S. Trofimchuk},
On the geometry of wave solutions of a delayed reaction-diffusion
equation, {E-print: arXiv:math/0611753v2 [math.DS] (2008) 25 p.}

\bibitem{ETST} E. Trofimchuk, S. Trofimchuk, Admissible wavefront speeds for a single species reaction-diffusion equation with delay, {
Discrete Contin. Dynam. Systems A,} { 20} (2008) 407-423.

\bibitem{WLR} Z.C. Wang, W.T. Li, S. Ruan, Traveling fronts in monostable equation with nonlocal delayed effects , {
J. Dyn. Diff. Equat.,} published online: 11 March 2008.

\bibitem{WM} J. Wu, D. Wei, M. Mei, Analysis on the critical speed of traveling waves, {
Applied Mathematics Letters,} { 20} (2007) 712-718.



\end{thebibliography}
\end{document}